\documentclass{amsart}
\usepackage{amssymb}
\usepackage{amsbsy}
\usepackage{amscd}
\usepackage{amsmath}
\usepackage{amsthm}
\usepackage{amsxtra}
\usepackage{latexsym}
\usepackage[mathscr]{eucal}
\usepackage{upgreek}
\usepackage{wasysym}
\usepackage[all,cmtip]{xy}

\usepackage{xcolor}
\definecolor{rouge}{rgb}{0.85,0.1,.4}
\definecolor{bleu}{rgb}{0.1,0.2,0.9}
\definecolor{violet}{rgb}{0.7,0,0.8}
\usepackage[dvipdfmx,colorlinks=true,linkcolor=bleu,urlcolor=violet,citecolor=rouge]{hyperref}

\newcommand{\bra}{{\langle}}
\newcommand{\ket}{{\rangle}}

\newcommand{\cprime}{$'$}
\newcommand{\on}{\operatorname}
\newcommand{\+}{\mathop{\oplus}}
\renewcommand{\*}{{\otimes}}
\newcommand{\mc}{\mathcal}
\newcommand{\mf}{\mathfrak}

\newcommand{\g}{\mf{g}}

\newcommand{\affg}{\widehat{\mf{g}}}

\newcommand{\Z}{\mathbb{Z}}
\newcommand{\C}{\mathbb{C}}
\newcommand{\N}{\mathbb{N}}

\newcommand{\W}{\mathscr{W}}
\newcommand{\ra}{\rightarrow}
\newcommand{\lam}{\lambda}

\def\geq{\geqslant}

\DeclareMathOperator{\gr}{gr}

\theoremstyle{theorem}
\newtheorem{Th}{Theorem}[section]

\newtheorem{Lem}[Th]{Lemma}

\theoremstyle{remark}

\newtheorem{Conj}{Conjecture}

\title{Associated varieties and Higgs branches (A survey)}
\author{Tomoyuki Arakawa}
\address{Research Institute for Mathematical Sciences, Kyoto University,
Kyoto 606-8502 JAPAN}
\address{Department of Mathematics,
Massachusetts Institute of Technology,
77 Massachusetts Avenue
Cambridge, MA 02139-4307 USA}
\email{arakawa@kurims.kyoto-u.ac.jp}

\begin{document}
\begin{abstract}
Associated varieties of vertex algebras are 
analogue of the associated varieties of primitive ideals of the universal enveloping algebras of 
semisimple Lie algebras.
They not only capture  some of the important properties of vertex algebras 
but also have interesting relationship with the Higgs branches of 
four-dimensional $N=2$ superconformal field theories (SCFTs).
As a consequence, one can deduce the modular invariance of Schur indices 
of 4d $N=2$ SCFTs
from  the theory of vertex algebras.\end{abstract}
\maketitle

%

\section{Associated varieties of vertex algebras}
A \emph{vertex algebra} consists of a vector space $V$ with a distinguished vacuum vector $|0\ket  \in V$ and a vertex operation, which is a linear map $V \otimes V \rightarrow V((z))$, written $u \otimes v \mapsto Y(u,z)v = (\sum_{n \in \Z} u_{(n)} z^{-n-1})v$, such that the following are satisfied:
\begin{itemize}
\item (Unit axioms) $(|0\ket)(z) = 1_V$ and $Y(u,z)|0\ket  \in u + zV[[z]]$ for all $u \in V$.

\item (Locality)
$(z-w)^n[Y(u,z),Y(v,w)]=0$ for a sufficiently large $n$
for all $u, v\in V$.
\end{itemize}
The operator $\partial: u \mapsto u_{(-2)}|0\ket $ is called the translation operator and it satisfies $Y(Tu, z) = \partial_z Y(u, z)$. The operators $u_{(n)}$ are called \emph{modes}.

To each vertex algebra $V$
one associates a Poisson algebra
$R_V$,
called the {\em Zhu's $C_2$-algebra},
as follows (\cite{Zhu96}).
 Let $C_2(V)$ be the subspace of $V$ spanned by the elements
 $a_{(-2)}b$, $a,b\in V$,
 and 
set 
 $R_V=V/C_2(V)$.
Then $R_V$ is a Poisson algebra
by 
\begin{align*}
\bar a.\bar b=\overline{a_{(-1)}b},\quad \{ \bar a,\bar b\}=\overline{a_{(0)}b},
\end{align*}
where $\bar a$ denote the image of $a\in V$ in $R_V$.

A vertex algebra is called  {\em strongly finitely generated} if $R_V$ is finitely generated.
In this note we assume that all the vertex algebras are finitely strongly generated.

The 
{\em associated variety} $X_V$ of a vertex algebra $V$
is  the affine Poisson variety $X_V$ defined by
\begin{align*}
X_V=\on{Specm}(R_V)
\end{align*}
(\cite{Ara12}).

\smallskip

Let $\g$ be a simple Lie algebra over $\C$,
 $\affg=\g[t,t^{-1}]\+ \C K$ be the affine Kac-Moody algebra associated with $\g$
and the normalized invariant inner product $(~|~)$.
Set
$$V^k(\g):=U(\affg)\*_{U(\g[t]\+ \C K)}\C_k,$$
where $k\in \C$ and $\C_k$ is one-dimensional representation of $\g[t]\+ \C K$
on which $\g[t]$ acts trivially and $K$  acts as the multiplication by $k$.
There is a unique vertex algebra structure on
$V^k(\g)$ 
such that $|0\ket =1\*1$ is the vacuum vector
and $$Y(x,z)=x(z):=\sum_{n\in \Z}(xt^n)z^{-n-1}\quad (x\in \g),$$
where
we consider $\g$ as a subspace of $V^k(\g)$ by the embedding
$\g\hookrightarrow V^k(\g)$, $x\mapsto xt^{-1}|0\ket$.
$V^k(\g)$ 
 is called the {\em universal affine vertex algebra associated with $\g$
at level $k$}.

One can regard $V^k(\g)$ as an analogue of the universal enveloping algebra
in the sense that 
a $V^k(\g)$-module is the same as 
a smooth $\affg$-module of level $k$,
where a $\affg$-module $M$ is called {\em smooth} if 
$x(z)m\in M((z))$ for all $m\in M$, $x\in \g$,
and called  of level $k$ if $K$ acts as the multiplication by $k$.

Any graded quotient $V$ of $V^k(\g)$  as a $\affg$-module has the structure of the quotient vertex algebra.
In particular the unique simple graded quotient $L_k(\g)$ is a vertex algebra
and is called the {\em simple affine vertex algebra associated with $\g$ at level $k$}.

For any quotient vertex algebra $V$ of $V^k(\g)$,
we have
$R_V=V/\g[t^{-1}]t^{-2}V$,
and the surjective linear map
\begin{align}
\C[\g^*]=S(\g)\ra R_V,\quad  x_1\dots x_r\mapsto \overline{(x_1t^{-1}\dots x_rt^{-1}|0\ket}\quad (x_i\in \g)
\label{eq:Poisson-surj}
\end{align}
is a homomorphism of Poisson algebras.
In particular  $X_{L_k(\g)}$ is a subvariety of $\g^*$,
which is $G$-invariant and conic.
We note that on the contrary to the associated variety of a primitive ideal of $U(\g)$,
$X_{L_k(\g)}$  is not necessarily contained in the nilpotent cone $\mc{N}$ of $\g$.
Indeed,
$L_k(\g)=V^k(\g)$ for a generic $k$ and 
$X_{V^k(\g)}=\g^*$ as \eqref{eq:Poisson-surj} is an isomorphism
for $V=V^k(\g)$ by the PBW theorem.
\smallskip

For a nilpotent element $f$ of $\g$,
let $\W^k(\g,f)$ be the universal $W$-algebra associated with $(\g,f)$ at level $k$:
\begin{align*}
\W^k(\g,f)=H_{DS,f}^0(V^k(\g)),
\end{align*}
where $H_{DS,f}^\bullet(?)$ is the BRST cohomology functor of the quantized Drinfeld-Sokolov reduction
associated with $(\g,f)$ (\cite{FF90,KacRoaWak03}).
The associated variety 
$X_{\W^k(\g,f)}$ is isomorphic to the {\em Slodowy slice} $\mc{S}_f=f+\g^e$,
where $\{e,f,h\}$ is an $\mf{sl}_2$-triple and $\g^e$ is the centralizer of $e$ in $\g$ (\cite{De-Kac06}).
For any quotient $V$ of $V^k(\g)$,
$H_{DS,f}^0(V)$ is a quotient vertex algebra of $\W^k(\g,f)$ provided that it is nonzero,
and we have
\begin{align}\label{eq:ass-variety-W}
X_{H_{DS,f}^0(V)}= X_V\cap \mc{S}_f,
\end{align}
which is a $\C^*$-invariant subvariety of $\mc{S}_f$ (\cite{Ara09b}).

\section{Lisse and quasi-lisse vertex algebras}
A vertex algebra $V$ is called {\em lisse} (or $C_2$-cofinite) if 
$\dim X_V=0$,
or equivalently, $R_V$ is finite-dimensional.
For instance,  $L_k(\g)$ is  lisse if and only if $L_k(\g)$ is integrable as a $\affg$-module, 
or equivalently,
$k\in \Z_{\geq 0}$ (\cite{DonMas06}).
Therefore, the lisse condition generalizes the integrability to an arbitrary vertex algebra.
Indeed, lisse vertex algebras are analogue of finite-dimensional algebras in the following sense.
\begin{Lem}[\cite{Ara12}]
A vertex algebra $V$ is lisse if and only if $\dim \on{Spec}(\gr V)=0$,
where $\gr V$ is the associated graded Poisson vertex algebra with respect to the canonical filtration on $V$ (\cite{Li05}).
\end{Lem}
It is known that lisse vertex algebras have various nice properties
such as
modular invariance of characters of $V$-modules under some mild assumptions (\cite{Zhu96,Miy04}).
However, there are significant vertex algebras that do not satisfy the lisse condition.
For instance, an {\em admissible affine vertex algebra} $L_k(\g)$ 
(see below)
has a 
complete reducibility property (\cite{A12-2}) 
and the modular invariance property
(\cite{KacWak89}, see also \cite{AEkeren})
in the category $\mc{O}$ 
although it is not lisse
unless it is integrable.
So it is natural to try to relax the lisse condition.

Since
$X_V$ is a Poisson variety
we have a finite partition
$$X_V=\bigsqcup_{k=0}^r X_k,$$
where each $X_k$ is a smooth analytic Poisson variety.
Thus for any point $x\in X_k$ there is a well defined symplectic leaf through it.
A vertex algebra $V$ is called {\em quasi-lisse} (\cite{Arakawam:kq}) if $X_V$ has only finitely many symplectic leaves.
Clearly,
lisse vertex algebras are quasi-lisse.

For example, 
consider the simple affine vertex algebra $L_k(\g)$.
Since
symplectic leaves in $X_{L_k(\g)}$ are the coadjoint  $G$-orbits 
contained in $X_{L_k(\g)}$,
where $G$ is the adjoint group of $\g$,
 it follows that 
$L_k(\g)$ is quasi-lisse if and only if $X_{L_k(\g)}\subset \mc{N}$.
Hence \cite{FeiMal97,Ara09b},
admissible affine vertex algebras  are quasi-lisse.

A theorem of Etingof and Schelder \cite{EtiSch10}
says that if a Poisson variety $\on{Specm}(R)$ has finitely many symplectic leaves then
the zeroth Poisson homology
$R/\{R,R\}$ is finite-dimensional.
It follows \cite{Arakawam:kq} that 
a quasi-lisse conformal vertex algebra has only finitely many 
simple ordinary  representations.
Here 
a $V$-module $M$ is 
called ordinary if it is a positive energy representation
on which $L_0$ acts semisimply
and
 each $L_0$-eigenspace is finite-dimensional,
so that 
the normalized character 
\begin{align*}
\chi_M(\tau)=\on{tr}_M (q^{L_0-\frac{c}{24}})
\end{align*}
is well-defined.

By extending Zhu's argument \cite{Zhu96}
using the theorem of Etingof and Schelder,
we get the following assertion.
\begin{Th}[\cite{Arakawam:kq}]\label{Th:AK}
Let $V$ be a quasi-lisse vertex algebra and $M$ a ordinary $V$-module.
Then $\chi_M$  satisfies a {\em modular linear differential equation}.
\end{Th}
Since the space of solutions of a modular linear differential equation (MLDE) is invariant under the action of $SL_2(\Z)$,
this implies that a quasi-lisse vertex algebra possesses a certain  modular invariance property,
although we do not claim that the normalized characters of $V$-modules span the space of the solutions.
 
 \section{Irreducibility conjecture and Examples of quassi-lisse vertex algebras}
 Let $\hat{\Delta}^{re}$ be the set 
 of real roots
 of $\affg$,
 $\hat{\Delta}^{re}_+$ the set of real positive roots.
 For a weight $\lam$ of $\affg$,
 let $\hat{\Delta}(\lam)=\{\alpha\in \hat{\Delta}^{re}\mid 
 \bra \lam+\rho,\alpha^{\vee}\ket\in \Z\}$,
 the integral roots system of $\lam$.
 An irreducible highest weight representation
 $L(\lam)$ of $\affg$ with highest weight $\lam$ is called admissible if $\lam$ is regular dominant,
that is,
$\bra \lam+\rho,\alpha^{\vee}\ket\not \in \{0,-1,-2,\dots,\}$ for all positive $\alpha\in \Delta_+$,
and $\mathbb{Q}\hat{\Delta}(\lam)=\mathbb{Q}\hat{\Delta}^{re}$ (\cite{KacWak89}).
The simple affine vertex algebra $L_k(\g)$  is called admissible if 
it is admissible as a $\affg$-module.
This condition is equivalent to that
\begin{align*}
k+h^{\vee}=\frac{p}{q},\quad p,q\in \N,\  (p,q)=1,\ p\geq \begin{cases}h^{\vee}&\text{if }(r^{\vee},q)=1,\\
h&\text{if }(r^{\vee},q)\ne 1,
\end{cases}
\end{align*}
where $h$, $h^{\vee}$,
and $r^{\vee}$ is the Coxeter number, the dual Coxeter number, and the lacing number of $\g$,
respectively (\cite{KacWak08}).

 As we have already mentioned above 
an admissible affine vertex algebra $L_k(\g)$ is quasi-lisse,
that is, $X_{L_k(\g)}\subset \mc{N}$.
In fact, the following assertion holds.
\begin{Th}[\cite{Ara09b}]
For an admissible affine vertex algebra
 $L_k(\g)$,
 $X_{L_k(\g)}$ is an irreducible variety contained in $\mc{N}$, that is,
 there exits a nilpotent orbit $\mathbb{O}$ such that 
 $X_{L_k(\g)}=\overline{\mathbb{O}}$.
\end{Th}
See \cite{Ara09b} for a concrete description of the orbit $\mathbb{O}$ that appears in the above theorem.

 For $\g=\mf{sl}_2$, it is not difficult to check that  $L_k(\g)$ is quasi-lisse
 if and only if $L_k(\g)$ is admissible
 for a non-critical\footnote{If $k$ is critical, that is, if $k=-h^{\vee}$,
then $X_{L_k(\g)}= \mc{N}$ by \cite{FeiFre92,EisFre01,FreGai04} for all simple Lie algebra $\g$.
}
 $k$, see \cite{Mal90,GorKac07}.
 However, there are non-admissible affine vertex algebras that are quasi-lisse for higher rank $\g$.
 
 Recall that the {\em Deligne exceptional series} \cite{De96} is the sequence of simple Lie algebras
 $$A_1\subset A_2\subset G_2\subset D_4\subset F_4\subset E_6\subset E_7\subset E_8.$$
 Let $\mathbb{O}_{min}$ be the unique non-trivial nilpotent orbit of $\g$.
\begin{Th}[\cite{AM15}]\label{Th:AMJoseph}
Let $\g$ be a simple Lie algebra that belongs to the Deligne exceptional series,
and let $k$ be a rational number of the form $k=-h^{\vee}/6-1+n$, $n\in \Z_{\geq 0}$, such that $k\not \in \Z_{\geq 0}$.
Then $$X_{L_k(\g)}=\overline{\mathbb{O}_{min}}.$$
\end{Th}
For types $A_1$, $A_2$, $G_2$, $D_4$, $F_4$,
the simple affine vertex algebra
$L_k(\g)$ appearing Theorem \ref{Th:AMJoseph}
is admissible, and hence, the statement is the special case of \cite{Ara09b}.
However, 
this is not the case for for types $D_4$, $F_4$, $E_6$, $E_7$, $E_8$
and
Theorem \ref{Th:AMJoseph}
gives examples of non-admissible quasi-lisse affine vertex algebras

Except for $\g=\mf{sl}_2$, the classification problem of quasi-lisse affine vertex algebras is  wide open.
(See \cite{AraMor16,AraMor16b} for more for more examples lisse affine vertex algebras.)

All the associated varieties are irreducible in the above examples
of quasi-lisse affine vertex algebras.
We conjecture that this is true in general:
\begin{Conj}[\cite{AraMor16}]\label{Conj1}
The associated variety of an quasi-lisse conical vertex algebra is irreducible.
\end{Conj}

\smallskip

Recall the description of associated variety of $W$-algebras given by \eqref{eq:ass-variety-W}.
This implies that
 if $L_k(\g)$ is quasi-lisse
and $f\in X_{L_k(\g)}$,
then
the $W$-algebra $H_{DS,f}^0(L_k(\g))$ is quasi-lisse as well,
and so is its simple quotient $\W_k(\g,f)$.
In this way we obtain 
a huge number of quasi-lisse $W$-algebras.
(See \cite{{AraMor16b}} for the irreducibility of the corresponding associated varieites.)
Moreover, if  $X_{L_k(\g)}=\overline{G.f}$,
then $X_{H_{DS,f}^0(L_k(\g))}=\{f\}$ by the transversality of $\mc{S}_f$ to $G$-orbits,
so that  $\W_k(\g,f)$ is in fact lisse.
Thus, Conjecture \ref{Conj1} in particular says that 
a quasi-lisse affine vertex algebra produces exactly one lisse simple $W$-algebra.

Lisse $W$-algebras thus obtained from admissible affine vertex algebras contain all the {\em exceptional $W$-algebras}
discovered by  Kac and Wakimoto \cite{KacWak08} (\cite{Ara09b}),
in particular, the {\em minimal series principal $W$-algebras} \cite{FKW92},
which are natural generalization of minimal series Virasoro vertex algebras \cite{BPZ84}.
The rationality of 
the minimal series principal $W$-algebras  has  been recently
recently proved  by the author (\cite{A2012Dec}).

\section{BL${}^2$PR${}^2$ correspondence and Higgs branch conjecture}
In \cite{BeeLemLie15},
Beem,  Lemos, Liendo, Peelaers, Rastelli,
and van Rees have constructed a remarkable 
map
$$\Phi: \{\text{4d $N=2$ SCFTs}\}\ra \{\text{vertex algebras}\},$$
such that,
among other things,
the character of the vertex algebra
$\Phi(\mc{T})$ coincides  with the {\em Schur index}
of the corresponding   4d $N=2$ SCFT $\mc{T}$,
which is an important invariant.

How do vertex algebras coming from 4d $N=2$ SCFTs look like?
According to \cite{BeeLemLie15}, we have
$$c_{2d}=-12c_{4d},$$
where $c_{4d}$ and $c_{2d}$ are central charges of the 4d $N=2$ SCFT
and the corresponding vertex algebra, respectively.
Since the central charge is positive for a unitary theory,
this implies that the vertex algebras obtained in this way are never unitizable.
In particular integrable affine vertex algebras never appear by this correspondence.

The main examples of vertex algebras
considered in  \cite{BeeLemLie15}
are affine vertex algebras $L_k(\g)$ of types 
$D_4$, $F_4$, $E_6$, $E_7$, $E_8$ at level $k=-h^{\vee}/6-1$,
which are non-rational,
non-admissible quasi-lisse affine vertex algebras appeared in Theorem \ref{Th:AMJoseph}.
One can find more examples in the literature, see e.g.\ \cite{Beem:2015yu,BN1,
CorSha16,
BN2,
Dan,SXY,BLN}.

Now, there is another important invariant of a 4d $N=2$ SCFT $\mc{T}$,
called the {\em Higgs branch},
which we denote by $Higgs_{\mc{T}}$.
The Higgs branch $Higgs_{\mc{T}}$ 
is an affine  algebraic variety
that has the hyperK\"{a}hler structure in its smooth part.
In particular,
$Higgs_{\mc{T}}$  is a (possibly singular)
 symplectic  variety.

Let $\mc{T}$ be one of the  4d $N=2$ SCFTs studied in  \cite{BeeLemLie15}
such that  that $\Phi(\mc{T})=L_k(\g)$ with 
$k=h^{\vee}/6-1$
for types 
$D_4$, $F_4$, $E_6$, $E_7$, $E_8$ as above. It is known that $Higgs_{\mc{T}}=\overline{\mathbb{O}_{min}}$ ,
which equals to $X_{L_k(\g)}$ by Theorem \ref{Th:AMJoseph}.
It is expected that this is not just a coincidence.
\begin{Conj}[Beem and Rastelli \cite{BeeRas}]\label{Conj:Beem and Rastell}
For a 4d $N=2$ SCFT $\mc{T}$, we have
\begin{align*}
Higgs_{\mc{T}}=X_{\Phi(\mc{T})}.
\end{align*}
\end{Conj}
So we are expected to  recover the Higgs branch of a 4d $N=2$ SCFT
from the corresponding vertex algebra, which is 
a purely algebraic object!
Note that the associated variety of a vertex algebra is only a Poisson variety in general.
Physical intuition expects that they are all quasi-lisse,
and so 
vertex algebras that come from 4d $N=2$ SCFTs via the map $\Phi$ form some special subclass of quasi-lisse vertex algebras.

We note that Conjecture \ref{Conj:Beem and Rastell} is a physical conjecture
since the Higgs branch is not a mathematically defined object at the moment.
 The Schur index is not a mathematically defined object  either.
However, in view of \cite{BeeLemLie15} and Conjecture \ref{Conj:Beem and Rastell},
one can try to define both Higgs branches and Schur indeces
of 4d $N=2$ SCFTs
using vertex algebras.
We note that there is a close relationship between Higgs branches of 4d $N=2$ SCFTs
and {\em Coulomb branches} of three-dimensional  $N=4$ gauge theories whose mathematical definition has been recently  given by 
Braverman, Finkelberg and Nakajima \cite{BFN16,BraFinNak17} (see \cite{A17a,A17b}).

In view of
Conjecture \ref{Conj:Beem and Rastell},
Theorem \ref{Th:AK} implies that 
the Schur index of a 4d $N=2$ SCFT satisfies a  MLDE,
which is something that has been conjectured in physics (\cite{BeeRas}).

\subsection*{Acknowledgments}
This note is based on the talks given by the author at 
 AMS Special Session ^^ ^^ Vertex Algebras and Geometry,"
at the University of Denver, October 2016,
 and at
 ^^ ^^ Exact operator algebras in superconformal field theories", 
 at Perimeter Institute for Theoretical Physics, Canada, December 2016.
 He thanks the organizers of these conferences
 and the Simons Collaboration on the Non-perturbative Bootstrap.
 He  benefited greatly from discussion with  Christopher Beem,
 Madalena Lemos,
 Anne Moreau,
 Hiraku Nakajima,
 Takahiro Nishinaka,
 Wolfger Peelaers,
 Leonardo Rastelli,
 Shu-Heng Shao,
 Yuji Tachikawa, and
 Dan Xie.
This research was supported in part 
JSPS KAKENHI Grant Numbers 17H01086, 17K18724,
and
by Perimeter Institute for Theoretical Physics. Research at Perimeter Institute is supported by the Government of Canada through Industry Canada and by the Province of Ontario through the Ministry of Economic Development \& Innovation.

\newcommand{\etalchar}[1]{$^{#1}$}

\bibliographystyle{alpha}
\bibliography{/Users/tomoyuki/Documents/Dropbox/bib/math}

\begin{thebibliography}{BPRvR]}

\bibitem[A1]{Ara12}
Tomoyuki Arakawa.
\newblock A remark on the {$C_2$} cofiniteness condition on vertex algebras.
\newblock {\em Math. Z.}, 270(1-2):559--575, 2012.

\bibitem[A2]{Ara09b}
Tomoyuki Arakawa.
\newblock Associated varieties of modules over {K}ac-{M}oody algebras and
  {$C_2$}-cofiniteness of {W}-algebras.
\newblock {\em Int. Math. Res. Not.}, 2015:11605--11666, 2015.

\bibitem[A3]{A2012Dec}
Tomoyuki Arakawa.
\newblock Rationality of {W}-algebras: principal nilpotent cases.
\newblock {\em Ann. Math.}, 182(2):565--694, 2015.

\bibitem[A4]{A12-2}
Tomoyuki Arakawa.
\newblock Rationality of admissible affine vertex algebras in the category
  {$\mathcal{O}$}.
\newblock {\em Duke Math. J.}, 165(1):67--93, 2016.

\bibitem[A5]{A17a}
Tomoyuki Arakawa.
\newblock Representation theory of W-algebras
and Higgs branch conjecture.
\newblock submitted to the Proceeding of ICM 2018.

\bibitem[A6]{A17b}
Tomoyuki Arakawa.
\newblock Chiral algebras of class $\mathcal{S}$ and symplectic varieties.
\newblock in preparation.


\bibitem[AK]{Arakawam:kq}
Tomoyuki Arakawam and Kazuya Kawasetsu.
\newblock Quasi-lisse vertex algebras and modular linear differential
  equations.
\newblock arXiv:1610.05865 [math.QA],
to appear in {\em Kostant Memorial Volume}, Birkhauser.

\bibitem[AM1]{AM15}
Tomoyuki Arakawa and Anne Moreau.
\newblock {J}oseph ideals and lisse minimal {W}-algebras.
\newblock {\em J. Inst. Math. Jussieu, {\em published online}.}

\bibitem[AM2]{AraMor16}
Tomoyuki Arakawa and Anne Moreau.
\newblock Sheets and associated varieties of affine vertex algebras.
\newblock {\em Adv.\ Math.} 
Vol.\ 320, 157--209.

\bibitem[AM3]{AraMor16b}
Tomoyuki Arakawa and Anne Moreau.
\newblock On the irreducibility of associated varieties of {W}-algebras.
\newblock to appear in 
 in the special issue of {\em J. Algebra} in Honor of Efim Zelmanov on occasion of his 60th anniversary.



\bibitem[AvE]{AEkeren}
Tomoyuki Arakawa and Jethro van Ekeren.
\newblock Modularity of relatively rational vertex algebras and fusion rules of
  regular affine {W}-algebras.
\newblock arXiv:1612.09100[math.RT].


\bibitem[BLL{\etalchar{+}}]{BeeLemLie15}
Christopher Beem, Madalena Lemos, Pedro Liendo, Wolfger Peelaers, Leonardo
  Rastelli, and Balt~C. van Rees.
\newblock Infinite chiral symmetry in four dimensions.
\newblock {\em Comm. Math. Phys.}, 336(3):1359--1433, 2015.

\bibitem[BPRvR]{Beem:2015yu}
Christopher Beem, Wolfger Peelaers, Leonardo Rastelli, and Balt~C. van Rees.
\newblock Chiral algebras of class {$\mathcal{S}$}.
\newblock {\em J. High Energy Phys.}, (5):020, front matter+67, 2015.

\bibitem[BR]{BeeRas}
Christopher Beem and Leonardo Rastelli.
\newblock Vertex operator algebras, {H}iggs branches, and modular differential
  equations.
\newblock arXiv:1707.07679[hep-th].

\bibitem[BPZ]{BPZ84}
A.~A. Belavin, A.~M. Polyakov, and A.~B. Zamolodchikov.
\newblock Infinite conformal symmetry in two-dimensional quantum field theory.
\newblock {\em Nuclear Phys. B}, 241(2):333--380, 1984.

\bibitem[BFN1]{BFN16}
Alexander Braverman, Michael Finkelberg, and Hiraku Nakajima.
\newblock Towards a mathematical definition of Coulomb branches of 3-dimensional $N=2$ gauge theories, II
\newblock 	arXiv:1601.03586 [math.RT].

\bibitem[BFN2]{BraFinNak17}
Alexander Braverman, Michael Finkelberg, and Hiraku Nakajima.
\newblock Ring objects in the equivariant derived satake category arising from
  coulomb branche.
\newblock {arXiv:1706.02112 math.RT]}.

\bibitem[BLN]{BLN}
Matthew Buican, Zoltan Laczko and Takahiro Nishinaka.
\newblock
$\mc{B} =2$ $\mf{S}$-duality revisited.
\newblock
{\em J. High Energ. Phys.} (2017) 2017: 87.

\bibitem[BN1]{BN1}
Matthew Buican and Takahiro Nishinaka.
\newblock On the superconformal index of Argyres--Douglas theories.
\newblock {\em J. Phys. A: Math. Theor.} 49 (2016), no. 1, 015401, 33 pp. 

\bibitem[BN2]{BN2}
Matthew Buican and Takahiro Nishinaka.
\newblock Conformal Manifolds in Four Dimensions and Chiral Algebras.
\newblock arXiv:1603.00887 [hep-th].


\bibitem[CS]{CorSha16}
Clay C{\'o}rdova and Shu-Heng Shao.
\newblock Schur indices, {BPS} particles, and {A}rgyres-{D}ouglas theories.
\newblock {\em J. High Energy Phys.}, (1):040, front matter+37, 2016.

\bibitem[Del]{De96}
Pierre Deligne.
\newblock La s\'erie exceptionnelle de groupes de {L}ie.
\newblock {\em C. R. Acad. Sci. Paris S\'er. I Math.}, 322(4):321--326, 1996.

\bibitem[DSK]{De-Kac06}
Alberto De~Sole and Victor~G. Kac.
\newblock Finite vs affine {$W$}-algebras.
\newblock {\em Japan. J. Math.}, 1(1):137--261, 2006.

\bibitem[DM]{DonMas06}
Chongying Dong and Geoffrey Mason.
\newblock Integrability of {$C\sb 2$}-cofinite vertex operator algebras.
\newblock {\em Int. Math. Res. Not.}, pages Art. ID 80468, 15, 2006.

\bibitem[EF]{EisFre01}
David Eisenbud and Edward Frenkel.
\newblock Appendix to \cite{Mus01}.
\newblock 2001.

\bibitem[ES]{EtiSch10}
Pavel Etingof and Travis Schedler.
\newblock Poisson traces and {$D$}-modules on {P}oisson varieties.
\newblock {\em Geom. Funct. Anal.}, 20(4):958--987, 2010.
\newblock With an appendix by Ivan Losev.

\bibitem[FF]{FF90}
Boris Feigin and Edward Frenkel.
\newblock Quantization of the {D}rinfel\cprime d-{S}okolov reduction.
\newblock {\em Phys. Lett. B}, 246(1-2):75--81, 1990.

\bibitem[FF]{FeiFre92}
Boris Feigin and Edward Frenkel.
\newblock Affine {K}ac-{M}oody algebras at the critical level and {G}el\cprime
  fand-{D}iki\u\i\ algebras.
\newblock In {\em Infinite analysis, Part A, B (Kyoto, 1991)}, volume~16 of
  {\em Adv. Ser. Math. Phys.}, pages 197--215. World Sci. Publ., River Edge,
  NJ, 1992.

\bibitem[FG]{FreGai04}
Edward Frenkel and Dennis Gaitsgory.
\newblock {$D$}-modules on the affine {G}rassmannian and representations of
  affine {K}ac-{M}oody algebras.
\newblock {\em Duke Math. J.}, 125(2):279--327, 2004.

\bibitem[FKW]{FKW92}
Edward Frenkel, Victor Kac, and Minoru Wakimoto.
\newblock Characters and fusion rules for {$W$}-algebras via quantized
  {D}rinfel\cprime d-{S}okolov reduction.
\newblock {\em Comm. Math. Phys.}, 147(2):295--328, 1992.

\bibitem[FM]{FeiMal97}
Boris Feigin and Fyodor Malikov.
\newblock Modular functor and representation theory of {$\widehat{\rm sl}_2$}
  at a rational level.
\newblock In {\em Operads: Proceedings of Renaissance Conferences (Hartford,
  CT/Luminy, 1995)}, volume 202 of {\em Contemp. Math.}, pages 357--405,
  Providence, RI, 1997. Amer. Math. Soc.

\bibitem[GK]{GorKac07}
Maria Gorelik and Victor Kac.
\newblock On simplicity of vacuum modules.
\newblock {\em Adv. Math.}, 211(2):621--677, 2007.

\bibitem[KRW]{KacRoaWak03}
Victor Kac, Shi-Shyr Roan, and Minoru Wakimoto.
\newblock Quantum reduction for affine superalgebras.
\newblock {\em Comm. Math. Phys.}, 241(2-3):307--342, 2003.

\bibitem[KW1]{KacWak89}
V.~G. Kac and M.~Wakimoto.
\newblock Classification of modular invariant representations of affine
  algebras.
\newblock In {\em Infinite-dimensional Lie algebras and groups
  (Luminy-Marseille, 1988)}, volume~7 of {\em Adv. Ser. Math. Phys.}, pages
  138--177. World Sci. Publ., Teaneck, NJ, 1989.

\bibitem[KW2]{KacWak08}
Victor~G. Kac and Minoru Wakimoto.
\newblock On rationality of {$W$}-algebras.
\newblock {\em Transform. Groups}, 13(3-4):671--713, 2008.

\bibitem[Li]{Li05}
Haisheng Li.
\newblock Abelianizing vertex algebras.
\newblock {\em Comm. Math. Phys.}, 259(2):391--411, 2005.

\bibitem[Mal]{Mal90}
F.~G. Malikov.
\newblock Verma modules over {K}ac-{M}oody algebras of rank {$2$}.
\newblock {\em Algebra i Analiz}, 2(2):65--84, 1990.

\bibitem[Miy]{Miy04}
Masahiko Miyamoto.
\newblock Modular invariance of vertex operator algebras satisfying {$C\sb
  2$}-cofiniteness.
\newblock {\em Duke Math. J.}, 122(1):51--91, 2004.

\bibitem[Mus]{Mus01}
Mircea Musta{\c{t}}{\u{a}}.
\newblock Jet schemes of locally complete intersection canonical singularities.
\newblock {\em Invent. Math.}, 145(3):397--424, 2001.
\newblock With an appendix by David Eisenbud and Edward Frenkel.


\bibitem[Ras]{Rastelli16}
Leonardo Rastelli.
\newblock Vertex operator algebras, higgs branches and modular differential
  equations.
\newblock In {\em String Math, Paris, June 28}, 2016.


\bibitem[SXY]{SXY}
Jaewon Song, Dan Xie, and Wenbin Yan.
\newblock Vertex operator algebras of Argyres-Douglas theories from M5-branes.
\newblock  {\em preprint}. 	
\newblock  arXiv:1706.01607[hep-th].

\bibitem[XYY]{Dan}
Dan Xie, Wenbin Yan, and Shing-Tung Yau.
\newblock Chiral algebra of {A}rgyres-{D}ouglas theory from {M}5 brane.
\newblock {\em preprint}.
\newblock arXiv:1604.02155[hep-th].

\bibitem[Zhu]{Zhu96}
Yongchang Zhu.
\newblock Modular invariance of characters of vertex operator algebras.
\newblock {\em J. Amer. Math. Soc.}, 9(1):237--302, 1996.

\end{thebibliography}

\end{document}